\documentclass{article}

\usepackage{amssymb}
\usepackage{amsmath}
\usepackage{graphicx}
\usepackage{amsthm}
\usepackage{indentfirst}
\usepackage{color}
\definecolor{small}{rgb}{0,0,1}
\definecolor{notsmall}{rgb}{0.05,0.75,0.2}
\input xy
\xyoption{all}

\newtheorem{defn}{Definition}[section]
\newtheorem{theo}[defn]{Theorem}
\newtheorem{prop}[defn]{Proposition}
\newtheorem{rem}[defn]{Remark}

\newtheorem{exa}[defn]{Example}
\newenvironment{dem}{\textbf{Proof.}\;}{\hspace{\stretch{1}}\rule{1ex}{1ex}}

\newcommand{\spa}{\mathbb{P}^3(\mathbb{C})}
\newcommand{\pla}{\mathbb{P}^2(\mathbb{C})}
\newcommand{\lin}{\mathbb{P}^1(\mathbb{C})}

\newcommand{\gen}{\mathbb{P}^{n-4}(\mathbb{C})}
\newcommand{\co}{\mathbb{C}}
\newcommand{\re}{\mathbb{R}}
\newcommand{\zn}{\mathbb{Z}}
\newcommand{\ssm}{X^{SS}(m)}
\newcommand{\sm}{X^S(m)}
\newcommand{\cq}{X^{SS}(m)//G}
\newcommand{\gq}{X^S(m)/G}
\newcommand{\gc}{\gamma^C}
\newcommand{\gl}{\gamma^L}

\newcommand{\Vl}{V^L}
\newcommand{\Vc}{V^C}
\newcommand{\mw}{\widehat{m}}
\newcommand{\cdi}{\mathrm{codim}}
\newcommand{\Gr}{\mathrm{Gr}}
\begin{document}

\title{Betti Numbers of GIT quotients of products of projective planes}

\author{ Francesca Incensi }

\date{November 2009}

\maketitle

\begin{abstract}
We study the GIT quotients for the diagonal action of the algebraic group $SL_3(\co)$
on the $n$-fold product of $\pla$: in particular we
determine a strategy in order to determine the (intersection) Poincar\'{e}
polynomial of any quotient variety. In the special case $n=6$ we determine an
explicit formula for the
(intersection) Betti numbers of a quotient variety,
depending only on the combinatorics of the weights of the polarization $m\in \zn^6_{>0}$.
\end{abstract}


\section*{Introduction}

\indent  Geometric Invariant Theory gives a method for constructing ``quotients''
varieties for any linear action of a complex reductive algebraic group
$G$ on a projective variety $X$ (see \cite{MuFoKi} and \cite{Dol} for the general setting).
The quotients depend on the choice of an ample linearized line bundle $L$:
in particular Dolgachev-Hu \cite{DoHu} and Thaddeus \cite{Tha} proved that only a
finite number of GIT quotients can be obtained when $L$ varies and
gave a general description of the maps relating the various
quotients.\\ In this paper we restrict our attention to the case when $G=\mathbb{P}SL_3(\co)$
acts diagonally on $X$, the $n$-fold product of $\pla$, $\pla^n$ (that we have already studied in \cite{Inc})
and describe the topology of an arbitrary quotient variety,
both geometric and categorical: as a main result we describe an algorithm to compute
the (intersection) Poincar\'{e} polynomial of any quotient variety. Moreover in the special case $n=6$
we obtain an explicit formula for the
(intersection) Betti numbers of a quotient variety,
depending only on the combinatorics of the weights of the polarization $m\in \zn^6_{>0}$.\\

The contents of the paper are more
precisely as follows.\\ \indent Section 1 is concerned with the main features of
the quotients of $G=\mathbb{P}SL_3(\mathbb{C})$ acting on $X=\pla^n$:
we remind the main features of a categorical quotient $\cq$ and we describe the $G$ ample cone $C^G(X)$.
At the end we discuss the birational maps $\theta_{\pm k}$ that can relate different quotients.\\
\indent Section 2 discusses small resolutions: we give conditions for a birational map $\theta_{+k}$ to be a small map,
and we discuss the problem of the existence of a small resolution.\\
\indent In Section 3 we give an algorithm that permits to compute the (intersecion) Poincar\'{e} polynomial of any quotient.
The main key is the decomposition theorem of Beilinson-Bernstein-Deligne.\\
\indent Section 4 contains the special case $n=6$: Theorem \ref{betti formula} contains an explicit formula for the
(intersection) Betti numbers of any quotient: this formula depends only on the combinatorics of the polarization $m\in \zn^6_{>0}$.


\section{$G=\mathbb{P}SL_3(\co)$ acting on
$X=\pla^n$}
First of all
remind that $\textrm{Pic}^G(X)\cong \zn^n\,$: an ample line bundle $L$
over $X$ is determined by $L=L(m):=L(m_1,\ldots,
m_n)$, where  $m_i\in\zn_{>0}\,,\;\forall i$ (supp. $m_i\geq m_{i+1}$ for all $i:1,\ldots,n-1$).
Then the set of
semi-stable points $\ssm$ is described by the Hilbert-Mumford numerical
criterion;
let $x=(x_1,\ldots,x_n)\in X$ and $|m|:=\sum_{i =1}^n m_i\,$:
\begin{equation}\label{eq:sss} x\in \ssm \Leftrightarrow
\left\{
\begin{array}{l}\sum_{k, x_k = y}m_k \leq \frac{|m|}{3}\\\\
\sum_{j, x_j \in r} m_j \leq 2\frac{|m|}{3}
\end{array}\right.\end{equation}
for every point $y\in\pla$ and for every line $r\subset\pla$.
Moreover $x$ is a stable point, $x \in \sm$, iff the numerical criterion (\ref{eq:sss})
is verified with strict inequalities.\\
\indent The numerical criterion can be restated as follows: if
$K$ and $J$ are subsets of
 $\{1,\ldots,n\}$, then we can associate them with
the numbers:
\[\gc_K(m)\,= |m|-3 \sum_{k \in K} m_k \,,\qquad
\gl_J(m)\, = 2|m|-3 \sum_{j \in J} m_j \,. \] In particular we have:
$ \gc_J(m)\,=\, -\gl_{J'}(m)$\; where $ J'=\{1,\ldots,n\}\setminus J$.\\
\indent
Now consider $K\subset\{1,\ldots,n\}$ with $|K|\geq2$ and the set $\Vc_K$
of configurations $(x_1,\ldots,x_n)\in X$ where the points indexed by $K$
are coincident and there are no further coincidences and no non-implied collinearities; in the same way,
if $J\subset\{1,\ldots,n\}, |J|>2$ consider the set $\Vl_J$
of configurations $(x_1,\ldots,x_n)\in X$ where the points indexed by $J$
are collinear and there are no further collinearities.\\
In order to study $\ssm$ and $\sm$, it is sufficient to consider the subsets $\Vc_K$ and $\Vl_J$: in fact
\[\Vc_K\subseteq \ssm\;\Leftrightarrow\; m_i\leq\frac{|m|}{3} \textrm{ for all } i\,
\textrm{ and } \gc_K(m)\geq 0\,,\]
and
\[\Vl_J\subseteq \ssm\;\Leftrightarrow\; m_i\leq\frac{|m|}{3} \textrm{ for all } i\,
\textrm{ and } \gl_J(m)\geq 0\,,\]
with similar statements for $\sm$ (for more details, see \cite{Inc}).\\
Moreover we observe that if $m_i<|m|/3$ for all $i$,
then the number of coincident points is at most $n-3$.\\
\indent
Now consider a point \[\xi \in
Z(m):=\left(\cq\right)\setminus \left(\gq\right)\,;\] this is the image in
$\cq$ of different, strictly semi-stable orbits, that all have
in their closure a closed, minimal orbit $Gx$, for a certain
configuration $x$ that has $|K|$ coincident points, and the others
$n-|K|$ collinear; by the numerical criterion, we get $\gc_K(m)=0$
and $\gl_{K'}(m)=0$, where $K$ indicates the coincident points,
while $K'=\{1,\ldots\}\setminus K$ indicates the collinear ones.
In this case $\phi(\overline{\Vc_K}\cup\overline{\Vl_{K'}})$ has dimension $n-|K|-3$,
where $\phi$ is the projection map to the quotient: $\phi:\ssm\rightarrow\cq$.\\\\
\indent Some quotients are particularly easy to compute: let's consider two examples.
\begin{exa}\label{example:easy1} Consider
a polarization $m=(s,s,s,s,1,1,\ldots,1,1)$ such that
\[\left\{\begin{array}{l}
0<s<\frac{1}{3}(4s+n-4)\\\\
\frac{1}{3}(4s+n-4)<2s<\frac{2}{3}(4s+n-4)\\\\
3s>\frac{2}{3}(4s+n-4)
\end{array}\right.\quad\Rightarrow\quad s>2(n-4)\,.\]
In this case the quotient $\cq$ is isomorphic to the product of
$(n-4)$ copies of $\pla$.
\end{exa}
\begin{exa}\label{example:easy2}
Consider the polarization $m=(s,s,1,1,\ldots,1)$ such that
\[\left\{
\begin{array}{l}
0<s<\frac{1}{3}(2s+n-2)\\\\
\frac{1}{3}(2s+n-2)<2s<\frac{2}{3}(2s+n-2)\\\\
n-3<\frac{1}{3}(2s+n-2)
\end{array}
\right.\quad\Rightarrow\quad n-\frac{7}{2}<s<n-2\,.\]
In this case the quotient $\cq$ is isomorphic to the product of $\gen\times\gen$.
\end{exa}

\indent Dolgachev and Hu proved in \cite{DoHu} that varying the line bundle $L$ only a finite number
of different quotients can be obtained: the space that parametrizes these quotients is the \emph{$G$-ample cone}
$C^G(X)$, the convex cone in $NS^G(X)\otimes\re$ spanned by ample $G$-linearized line bundles $L$ with $X^{SS}_L\neq\emptyset$,
where $NS^G(X)$ is the (N\'{e}ron-Severi) group of $G$-line bundles modulo homological equivalence.
This cone is subdivided in walls and chambers: a polarization $L$ lies on a wall if and only if $X^S_L\subsetneq X^{SS}_L$.
A chamber is a connected component of the complement of the union of walls: all the polarizations in a chamber define the same
set of stable points.\\ In our case
\[C^G(X)=\{(m_1,\ldots,m_n;\lambda)\in\re^n\times\re_+\,: \sum_{i=1}^n m_i=3\lambda,0\leq m_i\leq \lambda, i=1,\ldots,n\}.\]
Moreover $\sm\subsetneq\ssm$ if and only if there exists a subset $K\subset\{1,\ldots,n\}$ such that
\[\sum_{i\in K}m_i=\frac{|m|}{3};\] this is equivalent to the condition that $L(m)$ belongs to the hyperplane
\[H_{C_K,L_{K'}}=\left\{(m;\frac{|m|}{3})\in C^G(X)\,: \sum_{i\in K}m_i=\frac{|m|}{3}, \sum_{j\in K'}m_j=\frac{2}{3}|m|, K\subset\{1,\ldots,n\}\right\}.\]
This is a codimension$-1$ wall.\\

\indent Now we want to introduce the birational maps that
may be constructed between two quotients:
they will be crucial in determining the Betti numbers of the quotients.\\
Let $m$ be a polarization such that $3$ divides $|m|$ and $\sm\neq
\emptyset, \sm \subsetneq \ssm$; let us consider ``variations'' of
$m$ as follows:
\[\mw = m \pm
(0,\ldots,0,\underbrace{1}_{k},0,\ldots,0)\,.\]
Then applying the numerical criterion we get $\ssm\subset X^S(\mw)=X^{SS}(\mw)\subset \ssm$:
these inclusions induce a morphism
\begin{equation}\label{theta} \theta_{\pm k}: X^S(\mw)/G \longrightarrow
\cq\,, \end{equation}
which is an isomorphism over $\gq$, while over
$Z(m)$ is a contraction of
subvarieties. Moreover if the quotient $\cq$ is singular (possible for $n\geq 6$),
then $\theta_{\pm k}$ is a resolution of singularities: from now on we assume $n\geq 6$.\\
\indent Consider a point $\xi \in Z(m)$: for what we have already observed, $\xi$
is determined by a closed, minimal orbit $Gx$ for a certain configuration $x$
that has $|K|$ coincident points and
the other $n-|K|$ collinear ($\gc_K(m)=\gl_{K'}(m)=0$).\\
\indent Now we want to calculate $d$, the dimension of $\theta_{\pm k}^{-1}(\xi)$:
by the numerical criterion, only
one between $\Vc_K$ and $\Vl_{K'}$ is included in $X^S(\mw)$.
\\\indent Dealing with an elementary transformation of ``plus''
type ($\mw \stackrel{+1_k}{\longrightarrow} m$), then
\begin{itemize}
\item[-]if $k \in K\, \Rightarrow
$
\begin{equation}\label{dim fiber i in J}
\theta^{-1}_{+k}(\xi)\cong\mathbb{P}^{n-|K|-3}(\co)\,,\quad d=n-|K|-3;
\end{equation}
\item[-] if $k \in K'\, \Rightarrow
$
\begin{equation}\label{dim fiber i in J'}
 \theta^{-1}_{+k}(\xi)\cong\mathbb{P}^{2|K|-3}(\co)\,,\quad d=2|K|-3\,.
\end{equation}
\end{itemize}
\indent Similarly, with an elementary transformation of ``minus'' type
($\mw \stackrel{-1_k}{\longrightarrow} m$), we have
\begin{equation}\label{dim fiber trasf2}
k \in K \Rightarrow  d=2|K|-3\,; \qquad
k \in K' \Rightarrow d=n-|K|-3\,.
\end{equation}
From now on we will consider only elementary transformation of ``plus'' type, $\theta_{+k}$;
in particular we want to study the properties of $\theta_{+k}$:
when $\theta_{+k}$ is a blow-up map? When a small map?

\section{Small resolutions}
\begin{defn} A proper surjective algebraic map $f:Y_1\rightarrow Y_2$
between irreducible complex $N$-dimensional algebraic varieties is \emph{small}
if $Y_1$ is nonsingular and, for all $r>0$,
\[\cdi_\co\{y\in Y_2\,|\, \dim_\co f^{-1}(y)\geq r\}>2r\,.\]
A \emph{small resolution} is a resolution of singularities which is a small map.
\end{defn}
Small maps are particularly relevant, because they preserve intersection homology.

\begin{prop}\label{small map proposition} Let $\cq$ be a quotient such that $m$ lies on a $1$ codimension wall
of $C^G(X)$ and consider a variation of the weights $\mw$; then the birational map
\[\theta_{+k}:X^S(\mw)/G\rightarrow \cq\,,\]
is a small resolution if for each $\xi\in Z(m)$, determined by a minimal closed orbit in
$\overline{\Vc_K}\cup\overline{\Vl_{K'}}$ ($2\leq |K|\leq n-3$),
\begin{equation}\label{small map conditions}
\dim(\theta_{+k}^{-1}(\xi))=\left\{
\begin{array}{lll}
2|K|-3\qquad& \textrm{if }\; 2\leq |K|< \frac{n+1}{3}\,,\\\\
n-|K|-3\qquad& \textrm{if }\;\frac{n-1}{3}<|K|\leq n-3\,.
\end{array}
\right.\end{equation}
\end{prop}
\begin{rem}
Studying the subdivision in chambers and walls of the $G$-ample cone
$C^G(X)$, the existence of a small resolution appears quite natural: consider a polarization $m$
that determines a singular quotient $\cq$ and lies on a codimension$-1$ wall:
each singularity of the quotient, determined by a certain $\phi(\overline{\Vc_K}\cup\overline{\Vl_{K'}}$),
can be solved by the birational morphisms
\[\theta_{+k}: X^S(m_+)/G\rightarrow\cq\quad \textrm{and}\quad \theta_{-k}:X^S(m_-)/G\rightarrow\cq\]
where the polarizations $m_+$ and $m_-$ lie on opposite chambers.
In order to have a small resolution we have to choose the ``right side'' of each wall.
\end{rem}
\begin{dem} The demonstration is based on the definition of \emph{small map}:
in order to verify the definition in our case, we have to consider the codimension of the sets
\[\left\{ \xi\in Z(m)\,|\;\dim(\theta_{+k}^{-1}(\xi))\geq r\right\}\] for $r>0$.
We have to distinguish two different cases:
\begin{enumerate}
\item if $\dim(\theta_{+k}^{-1}(\xi)=2|K|-3$, then
\[\begin{array}{rcl}\cdi\left\{\xi\in Z(m)\,|\, \dim(\theta_+^{-1}(\xi))\geq 2|K|-3 \right\}&>&2(2|K|-3)\\
2(n-4)-(n-|K|-3)&>&4|K|-6\\
3|K|&<&n+1\,.
\end{array}\]
\item if $\dim(\theta_{+k}^{-1}(\xi)=n-|K|-3$, then
\[\begin{array}{rcl}\cdi\left\{\xi\in Z(m)\,|\, \dim(\theta_{+k}^{-1}(\xi))\geq n-|K|-3 \right\}&>&2(n-|K|-3)\\
2(n-4)-(n-|K|-3)&>&2n-2|K|-6\\
3|K|&>&n-1\,.
\end{array}\]
\end{enumerate}
In other words $\theta_{+k}:X^S(\mw)\rightarrow\cq$ is a small map if for every $\xi \in Z(m)$, $\theta_{+k}^{-1}(\xi)$
has minimum dimension between $2|K|-3$ and $n-|K|-3$.
\end{dem}
\begin{rem}
Combining the previous result and formulas (\ref{dim fiber i in J}) and (\ref{dim fiber i in J'}),
we can restate Proposition \ref{small map proposition} as follows: in the same hypothesis of \ref{small map proposition},
the birational map $\theta_{+k}$
is a small resolution if for each $\xi\in Z(m)$, determined by a minimal closed orbit in
$\overline{\Vc_K}\cup\overline{\Vl_{K'}}$ ($2\leq |K|\leq n-3$),
\begin{equation}
\left\{\begin{array}{l}
k\in K'\\
2\leq |K|< \frac{n+1}{3}\,,\end{array}\right.\quad \textrm{or}\quad
\left\{\begin{array}{l}
k\in K\\
\frac{n-1}{3}<|K|\leq n-3\,.
\end{array}\right.
\end{equation}
\end{rem}
\begin{exa}\label{example:codim1}
Let us consider an example: $n=7, m=(9,4,4,4,4,4,1), |m|=30$;
this polarization lies on the codimension$-1$  wall $H_{C_{17},L_{23456}}$.\\
The map $\theta_{+6}:X^S(9,4,4,4,4,3,1)/G\rightarrow\cq$ is a small map,
while $\theta_{+1}:X^S(8,4,4,4,4,4,1)/G\rightarrow\cq$ is not a small map.
\end{exa} \indent As an immediate consequence we have $IH_\bullet(\cq)$ is isomorphic to
$H_\bullet(X^S(\mw))$ if $\theta_{+k}$ is a small map.\\
In general we have the following stronger result regarding the existence of small resolutions
whose proof is the same of \cite{Hu} Theorem 2.5, in the case of a torus action;
in fact using the Gelfand-MacPherson correspondence, the moduli space
$\left(\pla^{SS}\right)(m)//SL_3(\co)$ can be identified with a quotient in the Grassmannian
$\Gr(3,\co^n)$ acted on by the torus $\left(\co^*\right)^{n-1}$.
\begin{theo} For every singular quotient $\cq$, there exists a polarization $\widehat{m}$
such that $\widehat{m}$ lies on a chamber close enough to $m$ and
$\theta:X^S(\widehat{m})/G\rightarrow\cq$ is a small resolution.
\end{theo}
\begin{rem} We have already proved the result for $m$ lying on a codimesion$-1$ wall;
for higher codimension cases, the proof is based on the observation that
a codimension$-N$ wall is the intersection of $N$ codimension$-1$ walls $H_i$
and for each wall $H_i$ ($i=1,\ldots,N$)
$C^G(X)\setminus H_i$ has two connected components lying in two sides of $H_i$;
one of these components, named $C_{i,s}$, defines a small resolution according to the previous result.
Then the proof examines each wall, detects the right component $C_{i,s}$ and studies their intersection for $i=1,\ldots,N$.
\end{rem}
\begin{exa} Consider $n=8$ and the polarization $m=(19,16,7,7,6,2,2,1)$: it lies on a codimension$-3$ wall, given by
\[H_{C_{18},L_{234567}}\cap H_{C_{267},L_{13458}}\cap H_{C_{345},L_{12678}}\,.\]
A polarization $\widehat{m}$ such that $\theta:X^S(\widehat{m})/G\rightarrow\cq$ is a small map is
$\widehat{m}=( 19,16,7,6,6,2,1,1)$: in fact $\sm\subset X^S(\widehat{m})\subset \ssm$ and
$\theta^{-1}(\xi)$ has the correct dimension for each $\xi \in Z(m)$.
The polarization $\widehat{m}$ is obtained studying the three codimension$-1$
walls in order to detect the right components and their intersection.
\end{exa}

\section{Poincar\'{e} polynomial}
Let us consider a quotient $Y=Y(m)=\cq$ and define the \textsl{intersection Betti numbers} of $\cq$ as
\[ ib_i(Y):=\dim IH_i(Y,\mathbb{Q}), \quad i:0,\ldots,4(n-4)\,.\]
If $Y$ is non-singular, this definition coincides with the classical one.\\
In what follows we will use
\[P(Y)=\sum_i t^i \dim(H_i(Y))\] to denote the ordinary Poincar\'{e} polynomial of a variety $Y$ and use
\[IP(Y)=\sum_i t^i\cdot ib_i(Y)\] to denote the intersection Poincar\'{e} polynomial of $Y$.\\
\indent We have already stressed that small maps preserve intersection homology:
if $\theta_{\pm k}: X^S(\widehat{m})/G\rightarrow \cq$ is a small map, then
$H_\bullet(X^S(\mw)/G)\cong IH_\bullet(\cq),$ and as a consequence the Betti numbers and the Poincar\'{e} polynomial are preserved.\\
\indent In the previous section we have determined the conditions for a $\theta_{+k}$ map to be small, but
in the general case $\theta_{+k}$ is not small.
In this way in order to compute the Betti numbers for a general quotient
it is necessary to use the decomposition theorem of Beilinson-Bernstein-Deligne
(in particular we will consider the simplified version of \cite{Hu}):
\begin{theo} Let $f:X\rightarrow Y$ be a projective algebraic map and $X$ be a nonsingular variety. Then there exists
\begin{enumerate}
\item a stratification $Y=\bigcup_\alpha Y_\alpha$ of $Y$
\item a list of enriched strata $E_\beta=(Y_\beta, L_\beta)$
     where $Y_\beta$ is a stratum of $Y$ and $L_\beta$ is a local system over $Y_\beta$;
     moreover assume that evvery local system $L_\beta$ is trivial.
\end{enumerate}
Then there exists a collection of polynomials $\psi_\beta$ for all strata such that,
     \[P(X)=\sum_\beta IP(\overline{Y}_\beta)\cdot \psi_\beta\,,\]
 and for a point $y\in Y$
     \[IP(f^{-1}(y))=\sum_\beta IP_y(\overline{Y}_\beta)\cdot\psi_\beta\,.\]
\end{theo}
\indent Consider a polarization $m$ on a codimension$-1$ wall $H_{C_K,L_{K'}}$,
and $\widehat{m}$, $\overline{m}$ two polarizations close enough to $m$ lying in different sides of the wall $H_{C_K,L_{K'}}$:
\[\begin{array}{crclc}
X^S(\widehat{m})/G & & & & X^S(\overline{m})/G\\
                   &\searrow \widehat{\theta}& & \overline{\theta}\swarrow & \\
                   & & \cq & &
\end{array}
\]


\indent Consider the map $\widehat{\theta}:X^S(\widehat{m})/G\rightarrow\cq$: suppose that for $\xi \in Z(m)$,
$\widehat{\theta}^{-1}(\xi)\cong \mathbb{P}^{n-|K|-3}(\co)$.
Now apply the decomposition theorem to the map $\widehat{\theta}$, where the stratification of $\cq$ is given by
\[\cq= \gq\, \cup Z(m)\,.\]
Then there exist two polynomials $\psi_{0}$ and $\varphi_{1}$ such that
\[P(X^S(\widehat{m})/G)=\psi_{0}IP(\cq)+\psi_{1}IP(Z(m))\,.\]
These two polynomials may be determined:
\[\psi_0=1\,,\quad \psi_{1}=1+t^2+\ldots+t^{2(n-|K|-3)}-IP_\xi(\cq)\,,\]
where $\xi$ is a point in $Z(m)$. \\Then if $\xi$ is any point in $Z(m)$, we have that $P\left(X^S(\widehat{m})/G\right)$
is equal to
\begin{equation}\label{poincaré1}
IP(\cq)+\left(1+\ldots+t^{2(n-|K|-3)}-IP_\xi(\cq)\right)IP(Z(m)).
\end{equation}
\indent In the same way, study the map $\overline{\theta}:X^S(\overline{m})/G\rightarrow\cq$:
suppose that for $\xi \in Z(m)$,
$\overline{\theta}^{-1}(\xi)\cong \mathbb{P}^{2|K|-3}(\co)$.
\\Then if $\xi$ is any point in $Z(m)$, we have that $P\left(X^S(\overline{m})/G\right)$ is equal to
\begin{equation}\label{poincaré2}
IP(\cq)+\left(1+\ldots+t^{2(2|K|-3)}-IP_\xi(\cq)\right)IP(Z(m)).
\end{equation}
Subtracting (\ref{poincaré1}) and (\ref{poincaré2}), $P\left(X^S(\widehat{m})/G\right)-P\left(X^S(\overline{m})/G\right)$
is equal to
\[P\left(X^S(\widehat{m})/G\right)-P\left(X^S(\overline{m})/G\right)=
\varepsilon(H_{C_K,L_{K'}})Q_t(H_{C_K,L_{K'}})IP(Z(m))\,,\]
where $\varepsilon(H_{C_K,L_{K'}})$ and $Q_t(H_{C_K,L_{K'}})$ are defined by
\[
\varepsilon(H_{C_K,L_{K'}})=\left\{
\begin{array}{cll}
1       & \textrm{if }\quad n-|K|-3>2|K|-3\Rightarrow& |K|<\frac{n}{3}\\
-1\qquad  & \textrm{if }\quad n-|K|-3<2|K|-3\Rightarrow & |K|>\frac{n}{3}\\
0      & \textrm{if }\quad n-|K|-3=2|K|-3\Rightarrow & |K|=\frac{n}{3}
\end{array}\right.
\]
\[
Q_t(H_{C_K,L_{K'}})=\left\{
\begin{array}{cl}
t^{2(2|K|-3)+2}+\ldots+t^{2(n-|K|-3)}       & \textrm{if }\quad  |K|<\frac{n}{3}\\
t^{2(n-|K|-3)+2}+\ldots+t^{2(2|K|-3)}\quad  & \textrm{if }\quad  |K|>\frac{n}{3}\\
0      & \textrm{if }\quad |K|=\frac{n}{3}
\end{array}\right.
\]
We can also write
\begin{equation}\label{poincarèfinal}
P\left(X^S(\widehat{m})/G\right)=P\left(X^S(\overline{m})/G\right)+\varepsilon(H_{C_K,L_{K'}})Q_t(H_{C_K,L_{K'}})IP(Z(m))\,.
\end{equation}
The previous formula tells us how the Poincar\'{e} polynomial and the Betti numbers vary,
when we cross a codimension$-1$ wall $H_{C_K,L_{K'}}$;
in particular let study $IP(Z(m))$.\\
The locus $Z(m)$ is equal to $\cq\setminus\gq$ and it is determined by
\[\phi(\overline{\Vc_K}\cup\overline{\Vl_{K'}})\cong \left(\lin^{n-|K|}(m')\right)^{SS}//SL_2(\co)\,,\]
where $m'$ is obtained from $m$, by deleting all those weights $m_i$ with $i\in K$.
For the Hilbert-Mumford numerical criterion, the open set of stable points for $m'$
is equal to the open set of semi-stable points for $m'$, i.e. the categorical quotient
$\left(\lin^{n-|K|}(m')\right)^{SS}//SL_2(\co)$ is also geometric.\\
Now the Poincar\'{e} polynomials for this kind of quotient are well known (see \cite{HauKnu} for details):
\begin{equation}\label{poincarélin}
P(Z(m))=\frac{1}{1-t^2}\sum_{J\in S_{n-|K|}}\left(t^{2|J|}-t^{2(n-|K|-|J|-2)}\right)\,.
\end{equation}
where $S_{n-|K|}=\{J\subset\{1,2,\ldots,n-|K|\}\,:\, m'_{n-|K|}+\sum_{j\in J}m'_j < \sum_{i\notin J}m'_i\}$.
\begin{exa}In Example \ref{example:codim1} we have introduced the polarization
$m=(9,4,4,4,4,4,1)$ that
lies on the codimension$-1$  wall $H_{C_{17},L_{23456}}$.\\
We have $Z(m)\cong \left(\lin^5(4,4,4,4,4)\right)^{SS}//SL_2(\co)$ and
\[S_5=\{\emptyset, \{1\}, \{2\}, \{3\}, \{4\}\}.\]
\[\Rightarrow\; P(Z(m)=\frac{1}{1-t^2}(1-t^6+4t^2-4t^4)=1+5t^2+t^4\,.\]
In fact $Z(m)$ is isomorphic to the blow-up of $\pla$ in four points.
\end{exa}
\indent Now collecting all the results, we can calculate the Poincar\'{e} polynomial of any quotient
 $\cq$; this is the strategy:
\begin{enumerate}
\item[1)] If $m$ lies on a face of $C^G(X)$ and there exists one or more weights $m_i$ such that $m_i=0$,
           then the quotient has lower dimension than expected: study the polarization $m'$ obtained from $m$
           by deleting all the zero weights, and follow the next steps;
\item[2)] If $m$ lies on a face of $C^G(X)$ and there exists one weight $m_i$ such that $m_i=\frac{|m|}{3}$,
           then the quotient degenerates to $\left(\lin^{n-1}(m)'\right)^{SS}//SL_2(\co)$,
           where $m'$ is obtained from $m$ by deleting $m_i$; if $m'$ determines a categorical quotient
           that is also geometric use the Hausmann-Knutson formula, otherwise apply the decomposition theorem to this case;
\item[3)] If $m$ lies on a face of $C^G(X)$ and there exists two weights $m_i, m_j$ such that $m_i=m_j=\frac{|m|}{3}$,
           then the quotient degenerates to a point;
\item[4)] If $m$ lies on a wall (its codimension is not relevant),
           then determine a small resolution $X^S(\widehat{m})/G$ of $\cq$:
           the Poincar\'{e} polynomial does not change; now follow $5)$ for the polarization $\widehat{m}$;
\item[5)] If $m$ lies in a chamber, then construct a path $\gamma$ in $C^G(X)$
          that goes from a polarization $\widetilde{m}$ whose quotient is well-known (see Examples \ref{example:easy1}
          and \ref{example:easy2}) to $m$, and $\gamma$ meets only codimension$-1$ walls.
          Then if $\gamma$ crosses the codimension$-1$ walls $H_{C_{K_1},L_{K'_1}},\ldots, H_{C_{K_N},L_{K'_N}}$,
\[\gamma: m=\widetilde{m}=m_0 \stackrel{H_{C_{K_1},L_{K'_1}}}{\longrightarrow}  m_1 \stackrel{H_{C_{K_2},L_{K'_2}}}{\longrightarrow} m_2\rightarrow
 \ldots \stackrel{H_{C_{K_N},L_{K'_N}}}{\longrightarrow} m_N=m\]
          study each crossing and apply result (\ref{poincarèfinal}):
\[P\left(X^S(m_{j+1})/G\right)=P\left(X^S(m_j)/G\right)+\varepsilon(H_{C_{K_i},L_{K'_i}})Q_t(H_{C_{K_i},L_{K'_i}})IP(Z(m_{j,j+1})),\]
where $m_{j,j+1}$ indicates the polarization that lies on the wall $H_{C_{K_{j+1}},L_{K'_{j+1}}}$ and ``connects'' $m_j$ and $m_{j+1}$.
The Poincar\'{e} polynomial of each $Z(m_{j,j+1})$ can be computed using formula (\ref{poincarélin}).
\end{enumerate}
\section{A special case: $n=6$}
\indent Let's study the particular case $n=6$: it is really interesting because in this case every map
 $\theta_{+k}$ is the eventual composition of a small map and a blow-up map.\\
 \indent First of all remind that the number of chambers in which the $G$-ample cone $C^G(X)$
 is divided is less than or equal to 38 (see \cite{Inc} for details).\\


\indent For quotients $\cq$ such that $\sm\subsetneq \ssm$,
it is not so easy to give an upper bound to the number of quotients,
but we can get some important informations about their structure: the following result (\cite{Inc})
classifies the different types of points that may appear in $Z(m)$:
\begin{theo}\label{theo:sing} Let $X=\pla^6$ and $m\in
\zn_{>0}^6$ a polarization with $3\mid|m|$ and $m_i<|m|/3\, \forall i$: if
\begin{enumerate}
\item there are two different indexes $i,j$ s.t. $m_i+m_j=|m|/3$,
then the quotient includes a curve
$C_{ij}\cong \lin$, that corresponds to strictly semi-stable orbits
s.t. $x_i=x_j$ or $x_h, x_k, x_l, x_n$ collinear. In particular
points $\xi$ of $C_{ij}$ are singular: locally, the variety $(\cq,
\xi)$ is isomorphic to the toric variety\vspace{-0.2cm}
\[\co[T_1,T_2,T_3,T_4,T_5]/(T_1T_4-T_2T_3)\,.\]

\item there is
a ``partition'' of $m$ such that $m_i+m_j=m_h+m_l=m_k+m_n$,
then the quotient includes three curves
$C_{ij},$ $C_{hl},$ $C_{kn}\cong \lin$, that have a common point
$O_{ij,hl,kn}$. \\In particular $O_{ij,hl,kn}$ is singular: locally
the variety $(\cq, O_{ij,hl,kn})$ is isomorphic to the toric
variety\vspace{-0.2cm}
\[\co[T_1,T_2,T_3,T_4,T_5]/(T_1T_2T_3-T_4T_5)\,.\]
\item there are three indexes $h,i,j$ s.t. $m_h+m_i+m_j=|m|/3$,
then the quotient includes a point
$O_{hij}$ that correspond to the minimal, closed, strictly
semi-stable orbit $Gx$ such that $x_h=x_i=x_j$ and
$x_k, x_l, x_n$ are collinear. The point $O_{hij}$ is non singular.\\
\end{enumerate}
\end{theo}

\indent If $\xi$ is of type $3$ in Theorem \ref{theo:sing}, then $\theta_{\pm k}^{-1}(\xi)$
may be a point or $\spa$; in the last case $\theta_{\pm k}$ blows down $\spa$ to the point $\xi$:
\begin{equation}\label{retroimm theta}
\theta_{-k}^{-1}(\xi)\cong\left\{ \begin{array}{ll}\spa \quad & k\notin K\\
 \textrm{point} \quad & k\in K\end{array}\right.\qquad
 \theta_{+k}^{-1}(\xi)\cong\left\{ \begin{array}{ll}\textrm{point}  \quad & k\notin K\\
 \spa \quad & k\in K
\end{array}\right.
\end{equation}
\indent If $\xi$ is of type $1$, then $\xi \in C_{ij}\cong\lin$ and $\theta_{\pm k}^{-1}(C_{ij})$ has dimension two
(in particular it is one of the two dimensional quotients $\spa^5(m')/PSL_3(\co)$): it is a small contraction.\\
As the end we can consider $\theta_{\pm k}: X^S(\mw)/G \rightarrow
\cq$ as the eventual composition of a blow-down map with a small map; moreover if $\cq$ is singular, then
$\theta_{\pm k}$ is a resolution of singularities.\\

\indent The main result of this section is the following formula for the Poincar\'{e} polynomial of an arbitrary
categorical quotient $Y=\cq$:
\begin{theo} \label{betti formula}
Let $m=(m_1,\ldots,m_6)$ be a polarization such that
$0 < m_i < \frac{1}{3}\sum_i m_i\,, m_i\geq m_{i+1},$, then there may be five different cases:
\begin{enumerate}
\item if\; \[m_1+m_2+m_3 <\frac{2}{3}|m|,\] then $IP(Y)=1+6t^2+7t^4+6t^6+t^8$;
\item if \[\left\{\begin{array}{l}m_1+m_2+m_3 \geq\frac{2}{3}|m|\\\\ m_1+m_2+m_4<\frac{2}{3}|m|\end{array}\right.,\]
then $IP(Y)=1+5t^2+6t^4+5t^6+t^8$;
\item if \[\left\{\begin{array}{l}m_1+m_2+m_4 \geq\frac{2}{3}|m|\\\\ m_1+m_2+m_5<\frac{2}{3}|m|\\\\ m_1+m_3+m_4<\frac{2}{3}|m|\end{array}\right.,\]
 then $IP(Y)=1+4t^2+5t^4+4t^6+t^8$;
\item if \[\left\{\begin{array}{l}m_1+m_2+m_5 \geq\frac{2}{3}|m|\\\\ m_1+m_2+m_6<\frac{2}{3}|m|\end{array}\right. \quad or \quad
\left\{\begin{array}{l}m_1+m_3+m_4 \geq\frac{2}{3}|m|\\\\ m_2+m_3+m_4<\frac{2}{3}|m|\end{array}\right.,\]
then $IP(Y)=1+3t^2+4t^4+3t^6+t^8$;
\item if\; \[m_1+m_2+m_6\geq\frac{2}{3}|m|\quad or \quad m_2+m_3+m_4\geq\frac{2}{3}|m|,\]
then $IP(Y)=1+2t^2+3t^4+2t^6+t^8$.
\end{enumerate}
\end{theo}
\begin{dem} First of all let us observe that given two quotients $Y(m)$ and $Y(m')$
(geometric or categorical)
it is always possible to find
a finite sequence of algebraic maps $\theta_{\pm k}$ such that
$Y(m)\stackrel{\theta_{\pm i}}{\longrightarrow}\ldots
 \stackrel{\theta_{\pm j}}{\longrightarrow}Y(m')$, where each $\theta_{\pm k}$
 is a composition of a blow-down (or blow-up) map with a small map
 (or only one of these).\\
If $\theta_{+k}$ is a small map,
then it induces an isomorphism from the (intersection) homology
of $Y(\widehat{m})$ to the intersection homology of $Y(m)$.
It means that when we cross a wall of the $G$-ample cone, the intersection Betti numbers vary if the maps
$\theta_{\pm k}$ involve some blow-ups and blow-downs. \\As we have seen the map
$\theta_{\pm k}: Y(\mw) \longrightarrow Y(m)$ is a blow-down map if there is at least a point
$\xi=O_{hij}\in Y(m)$ of type $3$ (Theorem \ref{theo:sing}), that satisfies (\ref{retroimm theta}):
\[ \theta_{\pm k}^{-1}(O_{hij})=\theta_{\pm k}^{-1}\left(\phi(\overline{\Vc}_{hij}\cup\overline{\Vl}_{[6]\setminus\{h,i,j\}})\right)
=\widehat{\phi}(\overline{\Vl}_{[6]\setminus\{h,i,j\}})\cong\spa,\] where $\phi$ and $\widehat{\phi}$
are the projection to the quotients $Y(m)$ and $Y(\mw)$. In other words, $\theta_{\pm k}$ is a blow-down
if $X^S(\mw)$ contains $\Vl_{[6]\setminus\{h,i,j\}}$ and not $\Vc_{hij}$.\\
How many points $O_{hij}$ there may be in a categorical quotient $Y(m)$?
The answer is at most four and there are two possible sets of points:
\[\begin{array}{ll} O_{456},\quad O_{356},\quad O_{256},\quad O_{156},\qquad& \textrm{or}\\
O_{456},\quad O_{356},\quad O_{346},\quad O_{345}.&
\end{array}\]
For example if $m=(666621)$ then $O_{156}, O_{256}, O_{356}, O_{456}\in \cq$
(this quotient is particularly easy to compute: it is $\pla\times\pla$).\\
\indent In this way, in order to compute the intersection Betti numbers of a quotient $Y(m)$
(both geometric and categorical) it is sufficient to know the intersection Betti numbers
for a quotient $Y(m')$ and then check how many $\Vl_K$ sets (with $|K|=3$)
are NOT included in $\sm$: in fact for what we have just observed, we can always ``connect'' two quotients by
a finite sequence of maps that change the Betti numbers if and only if they are blow-ups or blow-down
of at most four copies of $\spa$. These copies of $\spa$ are detected by the sets $\Vl_K$ that are contained in $\sm$, where $K$ is one of the following:
\[456, 356, 256, 156\quad \textrm{or}\quad 456, 356, 346, 345\,.\]
\indent Consider now the chamber of $C^G(X)$ that contains the polarization $(5,5,5,5,1,1)$;
 then for every polarization $m$ in this chamber we have $Y(m)=\pla\times\pla$ and
moreover there are four sets $\Vl_K$ that are not included in $\sm$.
Its intersection Betti numbers are well-known and now we are able to compute all the Betti numbers:
\begin{itemize}
\item[-] if $m'$ is such that $X^S(m')$ does NOT contain four sets $\Vl_K$ then its intersection Betti numbers are the same of $\pla\times\pla$,
\item[-] if $m'$ is such that $X^S(m')$ does NOT contain three sets $\Vl_K$ then its intersection Betti numbers are the same of $\pla\times\pla$ blown-up in one point,
\item[-] if $m'$ is such that $X^S(m')$ does NOT contain two sets $\Vl_K$ then its intersection Betti numbers are the same of $\pla\times\pla$ blown-up in two points,
\item[-] if $m'$ is such that $X^S(m')$ does NOT contain one sets $\Vl_K$ then its intersection Betti numbers are the same of $\pla\times\pla$ blown-up in three points,
\item[-] if $m'$ is such that $X^S(m')$ contains all sets $\Vl_K$ then its intersection Betti numbers are the same of $\pla\times\pla$ blown-up in four points.
\end{itemize}
The relations of Theorem \ref{betti formula} compute the number of $\Vl_K$ sets that are included in $\sm$.
\end{dem}


\noindent \textsc{Francesca Incensi}\\
Dipartimento di Matematica, Universit\`{a} di Bologna, Italy\\
E-mail address: \textsl{incensi@dm.unibo.it}
\end{document}